\documentclass{article}
\usepackage{amssymb}
\usepackage{amsfonts}

\newtheorem{theorem}{Theorem}[section]

\newtheorem{corollary}{Corollary}

\newtheorem{definition}{Definition}[section]

\newtheorem{proposition}{Proposition}[section]
\newtheorem{remark}{Remark}[section]

\newenvironment{proof}[1][Proof]{\noindent\textbf{#1.} }{\ \rule{0.5em}{0.5em}}
\input{tcilatex}

\begin{document}

\title{AUTOMORPHIC EQUIVALENCE OF ONE-SORTED ALGEBRAS.\linebreak }
\author{{\Large A.Tsurkov} \\
%EndAName
\textit{Department of Mathematics and Statistics,} \\
\textit{Bar Ilan University,} \\
\textit{Ramat Gan, 52900, Israel.\bigskip }\\
\textit{Jerusalem College of Technology, }\\
\textit{21 Havaad Haleumi, }\\
\textit{Jerusalem, 91160, Israel.}\\
tsurkov@jct.ac.il}
\maketitle

\begin{abstract}
One of the central questions of universal algebraic geometry is: when two
algebras have the same algebraic geometry? There are various interpretations
of the sentence "Two algebras have the same algebraic geometry". One of
these is automorphic equivalence of algebras, which is discussed in this
paper, and the other interpretation is geometric equivalence of algebras. In
this paper we consider very wide and natural class of algebras: one sorted
algebras from IBN variety. The variety $\Theta $ is called an IBM variety if
two free algebras $W\left( X\right) ,W\left( Y\right) \in \Theta $ are
isomorphic if and only if the powers of sets $X$ and $Y$ coincide. In the
researching of the automorphic equivalence of algebras we must study the
group of automorphisms of the category $\Theta ^{0}$ of the all finitely
generated free algebras of $\Theta $ and the group of its automorphisms $%
\mathrm{Aut}\Theta ^{0}$. An automorphism $\Upsilon $ of the category $%
\mathfrak{K}$ is called inner if it is isomorphic to the identity
automorphism or, in other words, if for every $A\in \mathrm{Ob}\mathfrak{K}$
there exists $s_{A}^{\Upsilon }:A\rightarrow \Upsilon \left( A\right) $
isomorphism of these objects of the category $\mathfrak{K}$ and for every $%
\alpha \in \mathrm{Mor}_{\mathfrak{K}}\left( A,B\right) $ the diagram%
\[
\begin{array}{ccc}
A & \overrightarrow{s_{A}^{\Upsilon }} & \Upsilon \left( A\right) \\ 
\downarrow \alpha &  & \Upsilon \left( \alpha \right) \downarrow \\ 
B & \underrightarrow{s_{B}^{\Upsilon }} & \Upsilon \left( B\right)%
\end{array}%
\]%
is commutative. By \cite[Theorem 2]{PlotkinZhitAutCat}, if $\Theta $ is an
IBN variety of one-sorted algebras, then every automorphism $\Psi \in 
\mathrm{Aut}\Theta ^{0}$ can be decomposed: $\Psi =\Upsilon \Phi $, where $%
\Upsilon ,\Phi \in \mathrm{Aut}\Theta ^{0}$, $\Upsilon $ is an inner
automorphism of and $\Phi $ is a strongly stable one (see \textbf{Definition %
\ref{str_stab_aut}}). In this situation every strongly stable automorphism
defines the other algebraic structure on every algebra $H\in \Theta $ , such
that the algebra $H^{\ast }$ with this structure also belongs to our variety 
$\Theta $ and (\textbf{Theorem \ref{H_H*}}) even automorphically equivalent
to the algebra $H$, i.e., has the same algebraic geometry. From this we
conclude the necessary and sufficient conditions for two algebras to be
automorphically equivalent. We formulate these conditions by using the
notion of geometric equivalence of algebras. It means that we reduce
automorphic equivalence of algebras to the simpler notion of geometric
equivalence. This paper is a continuation of the research which was started
in \cite{PlotkinZhitAutCat}.
\end{abstract}

\section{Introduction.\label{intro}}

\setcounter{equation}{0}

We denote by $\Omega $ the signature of algebras of the variety $\Theta $.
Let \linebreak $X_{0}=\left\{ x_{1},x_{2},\ldots ,x_{n},\ldots \right\} $ be
a countable set of symbols, $\mathfrak{F}\left( X_{0}\right) $ - set of all
finite subset of $X_{0}$. We will consider the category $\Theta ^{0}$, which
objects are all free algebras $W\left( X\right) $ of the variety $\Theta $
generated by the finite subsets $X\in \mathfrak{F}\left( X_{0}\right) $.
Morphisms of the category $\Theta ^{0}$ are homomorphisms of these algebras.

In universal algebraic geometry we consider the "set of equations" $T\subset
B^{2}$ for some $B\in \mathrm{Ob}\Theta ^{0}$ and we "resolve" these
equations in the $\mathrm{Hom}\left( B,H\right) $ - "affine space over the
algebra $H\in \Theta $". We denote $T_{H}^{\prime }=\left\{ \mu \in \mathrm{%
Hom}\left( B,H\right) \mid T\subset \ker \mu \right\} $. This is the set of
all solutions of the set of equations $T$. For every set of "points" of
affine space $R\subset \mathrm{Hom}\left( B,H\right) $ we can consider a
congruence of equations defined by this set: $R_{H}^{\prime
}=\bigcap\limits_{\mu \in R}\ker \mu $. Also for every set $T\subset B^{2}$
we can consider its algebraic closer according the algebra $H$: $%
T_{H}^{\prime \prime }$. The set $T\subset B^{2}$ is called $H$-closed if $%
T=T_{H}^{\prime \prime }$. $H$-closed set is always a congruence. The
latices of the $H$-closed congruences in the algebra $B\in \mathrm{Ob}\Theta
^{0}$ we denote $Cl_{H}(B)$. We can consider the category of coordinate
algebras connected with the algebra $H\in \Theta $. This category we denote
by $C_{\Theta }\left( H\right) $. The objects of this category is quotients
algebras $W\left( X\right) /T$, where $X\in \mathfrak{F}\left( X_{0}\right) $%
, $T\in Cl_{H}(W\left( X\right) )$. The morphisms of this category is the
homomorphism of algebras. This category describes the algebraic geometry of
the algebra $H$. An answer to the our central question: are two algebras
have the same algebraic geometry - we can obtain by various comparisons of
these categories. All these definitions we can see, for example, in \cite%
{PlotkinVarCat}, \cite{PlotkinNotions} and \cite{PlotkinSame}.

\begin{definition}
\label{def_geom_equiv}\textbf{\hspace{-0.08in}. }\emph{(\cite{PlotkinVarCat})%
} \textit{Let }$H_{1}$\textit{\ and }$H_{2}$\textit{\ be algebras in }$%
\Theta $\textit{. The algebras }$H_{1}$\textit{\ and }$H_{2}$\textit{\ are
called geometrically equivalent if }$Cl_{H_{1}}(B)=Cl_{H_{2}}(B)$\textit{\
for every }$B\in \mathrm{Ob}\Theta ^{0}$\textit{.}
\end{definition}

Algebras $H_{1}$\ and $H_{2}$\ are geometrically equivalent if and only if
the categories $C_{\Theta }\left( H_{1}\right) $ and $C_{\Theta }\left(
H_{2}\right) $ coincide.

The notion of weak geometric equivalence of algebras defined in \cite[%
Section 4]{PlotkinVarCat} together with the closely connected notion of
similarity of algebras. In this paper we will give to this notion more
natural name: automorphic equivalence of algebras.

\begin{definition}
\label{def_aut_equiv}\textbf{\hspace{-0.08in}. }\textit{Let }$H_{1}$\textit{%
\ and }$H_{2}$\textit{\ be algebras in }$\Theta $\textit{.\ The algebras }$%
H_{1}$\textit{\ and }$H_{2}$\textit{\ are called automorphically equivalent
if there exists an automorphism }$\Phi :\Theta ^{0}\rightarrow \Theta ^{0}$%
\textit{\ and for every }$B\in \mathrm{Ob}\Theta ^{0}$\textit{\ there exists
a bijection}%
\[
\alpha (\Phi )_{B}:Cl_{H_{1}}(B)\rightarrow Cl_{H_{2}}(\Phi (B)) 
\]%
\textit{such that these bijections are coordinated with automorphism }$\Phi $%
\textit{\ in this sense: if }$B_{1},B_{2}\in \mathrm{Ob}\Theta ^{0}$\textit{%
, }$\mu _{1},\mu _{2}\in \mathrm{Hom}\left( B_{1},B_{2}\right) $\textit{, }$%
T\in Cl_{H_{1}}(B_{2})$\textit{\ and}%
\[
\tau \mu _{1}=\tau \mu _{2} 
\]%
\textit{\ then }%
\[
\widetilde{\tau }\Phi \left( \mu _{1}\right) =\widetilde{\tau }\Phi \left(
\mu _{2}\right) , 
\]%
\textit{where }$\tau :B_{2}\rightarrow B_{2}/T$\textit{, }$\widetilde{\tau }%
:\Phi \left( B_{2}\right) \rightarrow \Phi \left( B_{2}\right) /\alpha (\Phi
)_{B_{2}}\left( T\right) $\textit{\ are natural epimorphisms.}
\end{definition}

In \cite{PlotkinNotions} proved that if bijections $\left\{ \alpha (\Phi
)_{B}\mid B\in \mathrm{Ob}\Theta ^{0}\right\} $ are coordinated with
automorphism $\Phi $, then they defined uniquely by $\Phi $.

By the method of \cite{PlotkinVarCat} it can be proved that algebras $H_{1}$%
\ and $H_{2}$ are automorphically equivalent if and only if exists the pair $%
\left( \Phi ,\Psi \right) $ where $\Phi :\Theta ^{0}\rightarrow \Theta ^{0}$
is an automorphism,\textit{\ }$\Psi :C_{\Theta }\left( H_{1}\right)
\rightarrow $\textit{\ }$C_{\Theta }\left( H_{2}\right) $ is an isomorphism
and these three conditions

\begin{enumerate}
\item[A.] $\Psi \left( W\left( X\right) /Id\left( H_{1},X\right) \right)
=W\left( Y\right) /Id\left( H_{2},Y\right) $, where $\Phi \left( W\left(
X\right) \right) =W\left( Y\right) $,

\item[B.] $\Psi \left( W\left( X\right) /T\right) =W\left( Y\right) /%
\widetilde{T}$, where $T\in Cl_{H_{1}}(W\left( X\right) )$, $\widetilde{T}%
\in Cl_{H_{2}}(W\left( Y\right) )$,

\item[C.] natural epimorphism $\overline{\tau }:W\left( X\right) /Id\left(
H_{1},X\right) \rightarrow W\left( X\right) /T$ is transformed to the
natural epimorphism $\Psi \left( \overline{\tau }\right) :W\left( Y\right)
/Id\left( H_{2},Y\right) \rightarrow W\left( Y\right) /\widetilde{T}$
\end{enumerate}

are fulfilled ($Id\left( H,X\right) =\bigcap\limits_{\varphi \in \mathrm{Hom}%
\left( W\left( X\right) ,H\right) }\ker \varphi $ is a minimal $H$-closed
congruence in the $\left( W\left( X\right) \right) ^{2}$). It should be
remarked that if the pair $\left( \Phi ,\Psi \right) $, which fulfills
condition A. - C. exists, then the isomorphism $\Psi $ is defined uniquely
by $\Phi $.

The basic facts about automorphic equivalence are these:

\begin{enumerate}
\item If algebras $H_{1}$\ and $H_{2}$\ are geometrically equivalent then
they are\textit{\ }automorphically equivalent with the $\Phi =id_{\Theta
^{0}}$ and $\alpha (\Phi )_{B}$ is the identity mapping on $%
Cl_{H_{1}}(B)=Cl_{H_{2}}(B)$.

\item If an automorphism $\Phi :\Theta ^{0}\rightarrow \Theta ^{0}$ provides
the automorphic equivalence of algebras $H_{1}$ and $H_{2}$, then the
automorphism $\Phi ^{-1}:\Theta ^{0}\rightarrow \Theta ^{0}$ provides the
automorphic equivalence of algebras $H_{2}$ and $H_{1}$ with the $\alpha
(\Phi ^{-1})_{B}=\left( \alpha \left( \Phi \right) _{\Phi ^{-1}\left(
B\right) }\right) ^{-1}$.

\item If an automorphism $\Phi _{1}:\Theta ^{0}\rightarrow \Theta ^{0}$
provides the automorphic equivalence of algebras $H_{1}$ and $H_{2}$, and an
automorphism $\Phi _{2}:\Theta ^{0}\rightarrow \Theta ^{0}$ provides the
automorphic equivalence of algebras $H_{2}$ and $H_{3}$, then $\Phi _{2}\Phi
_{1}$ provides the automorphic equivalence of algebras $H_{1}$ and $H_{3}$
with the $\alpha (\Phi _{2}\Phi _{1})_{B}=\alpha \left( \Phi _{2}\right)
_{\Phi _{1}\left( B\right) }\alpha \left( \Phi _{1}\right) _{B}$.

\item[4.] If an automorphism $\Phi :\Theta ^{0}\rightarrow \Theta ^{0}$
which provides the automorphic equivalence of algebras $H_{1}$ and $H_{2}$
is an inner automorphism, then $H_{1}$ and $H_{2}$ are geometrically
equivalent.
\end{enumerate}

Facts 1., 2., 3. and 4. were formulated in \cite{PlotkinNotions} for the
similarity of algebras, but they can be very easily established for the
automorphic equivalence. From 1., 2. and 3. we conclude that the automorphic
equivalence is a reflexive, symmetric and transitive relation.

We will use this well known fact that conjugation of the inner automorphism
of some category $\mathfrak{K}$ by the arbitrary automorphism of this
category is also an inner automorphism.

\section{Derived and verbal operations.\label{operations}}

\setcounter{equation}{0}

\subsection{Verbal (polynomial) operations.\label{verbal_operations}}

Before the explanation of the notion of the verbal operation we will
introduce the short notation, which will be widely used in this paper. In
this notation $k$-tuple $\left( c_{1},\ldots ,c_{k}\right) \in C^{k}$ ($C$
is an arbitrary set) we denote by single letter $c$ and we will even allow
ourself to write $c\in C$ instead $c\in C^{k}$ and to write "homomorphism $%
\alpha :A\ni a\rightarrow b\in B$" instead "homomorphism $\alpha
:A\rightarrow B$, which transforms $a_{i}$ to the $b_{i}$, where $1\leq
i\leq k$".

For every word $w\left( x\right) \in W\left( X\right) $, where $X=\left\{
x_{1},\ldots ,x_{k}\right\} \subset X_{0}$ and every $H\in \Theta $ we can
define a $k$-ary operation $w_{H}^{\ast }\left( h\right) =w\left( h\right) $
(in full notation $h=\left( h_{1},\ldots ,h_{k}\right) \in H^{k}$, $x=\left(
x_{1},\ldots ,x_{k}\right) \in \left( W\left( X\right) \right) ^{k}$) or,
more formal, $w_{H}^{\ast }\left( h\right) =\gamma _{h}\left( w\right) $,
where $\gamma _{h}$ is a well defined homomorphism $\gamma _{h}:W\left(
X\right) \ni x\rightarrow h\in H$ (in full notation: homomorphism $\gamma
_{h}:W\left( X\right) \rightarrow H$, which transforms $x_{i}$ to the $h_{i}$
for $1\leq i\leq k$). This operation we call the \textbf{verbal operation}
induced on the algebra $H$ by the word $w\left( x\right) \in W\left(
X\right) $. If we will be very precise, we must say that "word $w\left(
x\right) \in W\left( X\right) $" is actually is a class $\left[ w\left(
x\right) \right] _{\Lambda _{X}}$ of words in the absolutely free algebra of
our signature $F\left( X\right) $ generated by set of symbols $X$, which are
congruent to the word $w\left( x\right) $ according the congruence $\Lambda
_{X}$ of the all identities of the variety $\Theta $ in $\left( F\left(
X\right) \right) ^{2}$. We define an operation $w_{H}^{\ast }$ on the
algebra $H$ of our variety $\Theta $. So the result of substitution $w\left(
h\right) $, or, by other words, the image $\gamma _{h}\left( w\right) $ does
not depend on what word from the class $\left[ w\left( x\right) \right]
_{\Lambda _{X}}$ we take. Because of that, we use the expression "word $%
w\left( x\right) \in W\left( X\right) $" and will avoid the redundant
punctuality, which will only impede the explanation. Also we must remark
that if the word $w\left( x\right) $ is generated by the set $X^{\prime
}\subset X$ such as $X^{\prime }\neq X$, then some variables in the
operation $w_{H}^{\ast }$ are fictive (the results of this operation do not
depend on them), but all our consideration is valid in this case.

\begin{remark}
\label{homomorphisms}\textbf{\hspace{-0.08in}. }\emph{{By \cite[1.8, Lemma 8]%
{GratzerUniveralg}, if $H_{1},H_{2}$ are algebras of the variety $\Theta $
and $\varphi $ a homomorphism from $H_{1}$ to $H_{2}$, then $\varphi \left(
w_{H_{1}}^{\ast }\left( h\right) \right) =w_{H_{2}}^{\ast }\left( \varphi
\left( h\right) \right) $ for every $h\in H_{1}$, i. e., $\varphi $ will be
a homomorphism from $H_{1}$ to $H_{2}$ as algebras with signature $\Omega
\cup \left\{ w_{H_{1}}^{\ast }\right\} $ and $\Omega \cup \left\{
w_{H_{2}}^{\ast }\right\} $ respectively.} }
\end{remark}

\subsection{Derived operations.\label{der_operations}}

We have another way to define additional algebraic operations on the
arbitrary algebra. Let us have an algebra $C$ with the signature $\Omega $
and a bijection $s:C\rightarrow C$. For every $\omega \in \Omega $ we can
define the derived operation $\widetilde{\omega }_{C}$ by this way: $%
\widetilde{\omega }_{C}\left( c\right) =s\left( \omega _{C}\left(
s^{-1}\left( c\right) \right) \right) $ for every $c\in C$ ($\omega _{C}$ is
a realization of the operation $\omega $ in the algebra $C$). By this
definition $s$ will be an isomorphism from the algebra $C$ with operations $%
\omega _{C}$ to the algebra $\widetilde{C}$, which has same domain and
operations $\widetilde{\omega }_{C}$ ($\omega \in \Omega $). Operations $%
\widetilde{\omega }_{C}$ ($\omega \in \Omega $) we call "derived operations
induced on the $C$ by the bijection $s$".

Now we will interweave the notions of derived and verbal operations. Let us
have a system of bijections $\left\{ s_{B}:B\rightarrow B\mid B\in \mathrm{Ob%
}\Theta ^{0}\right\} $ which fulfills these two conditions:

\begin{enumerate}
\item[B1)] for every homomorphism $\alpha :A\rightarrow B$ ($A,B\in \mathrm{%
Ob}\Theta ^{0}$) the mappings $s_{B}\alpha s_{A}^{-1}$ and $s_{B}^{-1}\alpha
s_{A}$ is also homomorphisms from $A$ to $B$ and

\item[B2)] $s_{B}\mid _{X}=id_{X}$ for every $B=W\left( X\right) \in \mathrm{%
Ob}\Theta ^{0}$.
\end{enumerate}

If arity of $\omega \in \Omega $ is $k$, we take $X_{\omega }=\left\{
x_{1},\ldots ,x_{k}\right\} \subset X_{0}$. $A_{\omega }=W\left( X_{\omega
}\right) $ - free algebra in $\Theta $. We have that $\omega \left( x\right)
\in A_{\omega }$ ($x=\left( x_{1},\ldots ,x_{k}\right) $) so there exists $%
w_{\omega }\left( x\right) \in A_{\omega }$ such that%
\begin{equation}
w_{\omega }\left( x\right) =s_{A_{\omega }}\left( \omega \left( x\right)
\right) .  \label{der_veb_opr}
\end{equation}

For every $H\in \Theta $ we denote $\omega _{H}^{\ast }$ the verbal
operation induced on the algebra $H$ by $w_{\omega }\left( x\right) $. $%
H^{\ast }$ will be the algebra, which has the same domain as the algebra $H$
and its operations are $\left\{ \omega _{H}^{\ast }\mid \omega \in \Omega
\right\} $. As it was proved in \cite[Theorem 3]{PlotkinZhitAutCat} 
\begin{equation}
\omega _{B}^{\ast }=\widetilde{\omega }_{B}  \label{zhito_eqv}
\end{equation}%
for every $B\in \mathrm{Ob}\Theta ^{0}$, where $\widetilde{\omega }_{B}$ is
a derived operations induced on the $B$ by the bijection $s_{B}$. But $%
\widetilde{\omega }_{H}$ is not defined for $H\in \Theta \smallsetminus 
\mathrm{Ob}\Theta ^{0}$.

By (\ref{zhito_eqv}) for every $B=W\left( X\right) \in \mathrm{Ob}\Theta
^{0} $ we have that $s_{B}:B\rightarrow B^{\ast }$ is an isomorphism. So the
system of words $\left\{ w_{\omega }\left( x\right) \in A_{\omega }=W\left(
X_{\omega }\right) \mid \omega \in \Omega \right\} $ fulfills these two
conditions:

\begin{enumerate}
\item[Op1)] $X_{\omega }=\left\{ x_{1},\ldots ,x_{k}\right\} $, where $k$ is
an arity of $\omega $, for every $\omega \in \Omega $;

\item[Op2)] for every $B=W\left( X\right) \in \mathrm{Ob}\Theta ^{0}$ there
exists an isomorphism $\sigma _{B}:B\rightarrow B^{\ast }$ (algebra $B^{\ast
}$ has same domain as the algebra $B$ and its operations $\omega _{B}^{\ast
} $ are induced by $w_{\omega }\left( x\right) $ for every $\omega \in
\Omega $) such as $\sigma _{B}\mid _{X}=id_{X}$
\end{enumerate}

with $\sigma _{B}=s_{B}$. This system of words we denote $\mathfrak{W}\left(
S\right) $.

\section{Systems of bijections, systems of words and strongly stable
automorphisms.\label{systems}}

\setcounter{equation}{0}

Now let us have a system of words $W=\left\{ w_{\omega }\left( x\right) \in
A_{\omega }=W\left( X_{\omega }\right) \mid \omega \in \Omega \right\} $
which fulfills conditions Op1) and Op2), then we can take a system of
bijections $\mathfrak{S}\left( W\right) =\left\{ \sigma _{B}:B\rightarrow
B\mid B\in \mathrm{Ob}\Theta ^{0}\right\} $.

By Op2) $B^{\ast }$ is a free algebra in the $\Theta $ with generators $%
\sigma _{B}\left( x\right) =x$ ($B=W\left( X\right) $, $x\in X$).

By Section \ref{operations}, we can induce the operation $\omega _{H}^{\ast
} $ by $w_{\omega }\left( x\right) $ ($\omega \in \Omega $) on every algebra 
$H\in \Theta $. As above $H^{\ast }$ the algebra which has the same domain
as the algebra $H$ and the operations $\left\{ \omega _{H}^{\ast }\mid
\omega \in \Omega \right\} $. By \textbf{Remark \ref{homomorphisms}}, if $%
\varphi :H_{1}\rightarrow H_{2}$ ($H_{1},H_{2}\in \Theta $) is a
homomorphism then $\varphi :H_{1}^{\ast }\rightarrow H_{2}^{\ast }$ is also
a homomorphism.

\begin{proposition}
\label{inthe_variety}\textbf{\hspace{-0.08in}. }$H^{\ast }\in \Theta $%
\textit{\ for every }$H\in \Theta $\textit{.}
\end{proposition}

\begin{proof}
There exists $B\in \mathrm{Ob}\Theta ^{0}$, such that $H$ is an epimorphic
image of $B$. So, $H^{\ast }$ is an epimorphic image of $B^{\ast }$. But $%
B^{\ast }\in \Theta $, hence $H^{\ast }\in \Theta $.
\end{proof}

We can consider the signature $\Omega ^{\ast }=\left\{ \omega ^{\ast }\mid
\omega \in \Omega \right\} $. Between signatures $\Omega ^{\ast }$ and $%
\Omega $ there is a symmetry:

\begin{proposition}
\label{invers}\textbf{\hspace{-0.08in}. }\textit{In every algebra }$H\in
\Theta $ \textit{operations }$\omega \in \Omega $\textit{\ are verbal
operations defined by the system of words }$U=\left\{ u_{\omega }^{\ast
}\left( x\right) \mid \omega \in \Omega \right\} $\textit{\ written by the
signature }$\Omega ^{\ast }$\textit{. The system of words }$U$\textit{\
fulfils conditions Op1) and Op2) with the system of isomorphism }$\sigma
_{B}^{-1}:B^{\ast }\rightarrow B$\textit{.}
\end{proposition}

\begin{proof}
$\omega \left( x\right) \in A_{\omega }=W\left( X_{\omega }\right) $, where $%
X_{\omega }=\left\{ x_{1},\ldots ,x_{k}\right\} $ and $k$ is the arity of $%
\omega $, so there exists $u\left( x\right) \in A_{\omega }$ such that $%
\omega \left( x\right) =\sigma _{A}\left( u_{\omega }\left( x\right) \right) 
$. We denote by $u_{\omega }^{\ast }\left( x\right) $ the word, which we
receive from the word $u_{\omega }\left( x\right) $ by replacement of all
operations $\omega \in \Omega $ by operations $\omega ^{\ast }\in \Omega
^{\ast }$ respectively. $\sigma _{A_{\omega }}$ is an isomorphism from the
algebra $A_{\omega }$ with operations $\omega \in \Omega $ to the algebra $%
A_{\omega }^{\ast }$ with operations $\omega ^{\ast }\in \Omega ^{\ast }$
and $\sigma _{A\omega }\left( x\right) =x$ for every $x\in X_{\omega }$, so $%
\sigma _{A_{\omega }}\left( u_{\omega }\left( x\right) \right) =u_{\omega
}^{\ast }\left( x\right) $. $\omega \left( h\right) =\gamma _{h}\left(
\omega \left( x\right) \right) =\gamma _{h}\left( u_{\omega }^{\ast }\left(
x\right) \right) $, where $\gamma _{h}:A_{\omega }\ni x\rightarrow h\in H$
is a homomorphism from $A_{\omega }$ to $H$ and from $A_{\omega }^{\ast }$
to $H^{\ast }$. $A_{\omega }^{\ast }$ with operations $\omega ^{\ast }\in
\Omega ^{\ast }$ is a free algebra in the variety $\Theta $ with the set of
generators $X_{\omega }$. So the operations $\omega \in \Omega $\ are verbal
operations. Condition Op1) is fulfilled by constructions of the words $%
u_{\omega }^{\ast }\left( x\right) $. Condition Op2) is obvious.
\end{proof}

And by \textbf{Remark \ref{homomorphisms}}, in which we change $H_{1},H_{2}$
to $H_{1}^{\ast },H_{2}^{\ast }$ correspondingly and vice versa we have

\begin{corollary}
\label{invers_homomorph}\textbf{\hspace{-0.08in}. }\textit{If }$%
H_{1},H_{2}\in \Theta $\textit{\ and }$\varphi :H_{1}^{\ast }\rightarrow
H_{2}^{\ast }$\textit{\ is a homomorphism, then }$\varphi :H_{1}\rightarrow
H_{2}$\textit{\ is also a homomorphism.}
\end{corollary}

\setcounter{corollary}{0}

And now we can prove

\begin{proposition}
\label{sys_biject}\textbf{\hspace{-0.08in}. }\textit{The system of
bijections }$\mathfrak{S}\left( W\right) =\left\{ \sigma _{B}:B\rightarrow
B\mid B\in \mathrm{Ob}\Theta ^{0}\right\} $\textit{\ fulfills conditions B1)
and B2).}
\end{proposition}

\begin{proof}
Actually, we must only prove condition B1). If $\alpha :A\rightarrow B$ ($%
A,B\in \mathrm{Ob}\Theta ^{0}$), then, because the operations in $A^{\ast }$
and $B^{\ast }$ are induced by the same words, by \textbf{Remark \ref%
{homomorphisms}} $\alpha $ is a homomorphism from $A^{\ast }$ to $B^{\ast }$%
, so $\sigma _{B}^{-1}\alpha \sigma _{A}$ will be a homomorphism from $A$ to 
$B$. $\sigma _{B}\alpha \sigma _{A}^{-1}$ is a homomorphism from $A^{\ast }$
to $B^{\ast }$. So, by a \textbf{Corollary} \ref{invers_homomorph} from 
\textbf{Proposition \ref{invers}}, it is a homomorphism from $A$ to $B$.
\end{proof}

\begin{proposition}
\label{operators}\textbf{\hspace{-0.08in}. }\textit{If }$S=\left\{
s_{B}:B\rightarrow B\mid B\in \mathrm{Ob}\Theta ^{0}\right\} $\textit{\ is a
system of bijections, which fulfills conditions B1) and B2), then }$%
\mathfrak{S}\left( \mathfrak{W}\left( S\right) \right) =S$\textit{. If }$%
W=\left\{ w_{\omega }\left( x\right) \mid \omega \in \Omega \right\} $%
\textit{\ is a system of words, which fulfills conditions Op1) and Op2),
then }$\mathfrak{W}\left( \mathfrak{S}\left( W\right) \right) =W$\textit{.}
\end{proposition}

\begin{proof}
Let $S=\left\{ s_{B}:B\rightarrow B\mid B\in \mathrm{Ob}\Theta ^{0}\right\} $%
. $\mathfrak{W}\left( S\right) $ be a system of words defined by the formula
(\ref{der_veb_opr}). By Section \ref{der_operations}, this system fulfills
conditions Op1) and Op2) with $\sigma _{B}=s_{B}$ for every $B\in \mathrm{Ob}%
\Theta ^{0}$. So $\mathfrak{S}\left( \mathfrak{W}\left( S\right) \right) =S$.

Let $W=\left\{ w_{\omega }\left( x\right) \mid \omega \in \Omega \right\} $
be a system of words, which fulfills conditions Op1) and Op2). $\mathfrak{S}%
\left( W\right) =\left\{ \sigma _{B}:B\rightarrow B\mid B\in \mathrm{Ob}%
\Theta ^{0}\right\} $, where isomorphisms $\sigma _{B}$ are remembered in
condition Op2). Let $\mathfrak{W}\left( \mathfrak{S}\left( W\right) \right)
=\{w_{\omega }^{\prime }\left( x\right) \mid \omega \in \Omega \}$. $%
\mathfrak{W}\left( \mathfrak{S}\left( W\right) \right) $ is defined by
formula (2.1), so $\sigma _{A_{\omega }}\left( \omega \left( x\right)
\right) =w_{\omega }^{\prime }\left( x\right) $. Also, by Op2), $\sigma
_{A_{\omega }}\left( \omega \left( x\right) \right) =\omega _{A_{\omega
}}^{\ast }\left( \sigma _{A}\left( x\right) \right) =\omega _{A_{\omega
}}^{\ast }\left( x\right) =w_{\omega }\left( x\right) $. Therefore, $%
w_{\omega }\left( x\right) =w_{\omega }^{\prime }\left( x\right) $ for every 
$\omega \in \Omega $ and $\mathfrak{W}\left( \mathfrak{S}\left( W\right)
\right) =W$.
\end{proof}

Now we introduce one of the central notion of our paper: strongly stable
automorphism of the category $\Theta ^{0}$. Automorphisms of this kind are
closely connected with the system of words, which fulfills conditions Op1)
and Op2) via system of bijections, which \ fulfills conditions B1) and B2).

\begin{definition}
\label{str_stab_aut}\textbf{\hspace{-0.08in}. }\textit{Automorphism of the
category }$\Theta ^{0}$\textit{\ are called \textbf{strongly stable} if it
fulfills these three conditions}

\begin{enumerate}
\item[A1)] $\Phi $\textit{\ preserves all objects of }$\Theta ^{0}$\textit{,}

\item[A2)] \textit{there exists a system of bijections }$\left\{ s_{B}^{\Phi
}:B\rightarrow B\mid B\in \mathrm{Ob}\Theta ^{0}\right\} $\textit{\ such
that }$\Phi $\textit{\ acts on the morphisms of }$\Theta ^{0}$\textit{\ by
this system, i. e., }%
\begin{equation}
\Phi \left( \alpha \right) =s_{B}^{\Phi }\alpha \left( s_{A}^{\Phi }\right)
^{-1}  \label{biject_action}
\end{equation}%
\textit{for every }$\alpha \in Hom\left( A,B\right) $\textit{\ (}$A,B\in
Ob\Theta ^{0}$\textit{) and}

\item[A3)] $s_{B}^{\Phi }\mid _{X}=id_{X}$\textit{\ for every }$B=W\left(
X\right) \in \mathrm{Ob}\Theta ^{0}$\textit{.}
\end{enumerate}
\end{definition}

Obviously, if we have a strongly stable automorphisms $\Phi $, then the
system of bijections $\left\{ s_{B}^{\Phi }:B\rightarrow B\mid B\in \mathrm{%
Ob}\Theta ^{0}\right\} $ fulfills conditions B1) and B2). We must remark,
that if a strongly stable automorphisms $\Phi $ of $\Theta ^{0}$ fulfills
conditions A2) and A3) with the system of bijections $\left\{ s_{B}^{\Phi
}:B\rightarrow B\mid B\in \mathrm{Ob}\Theta ^{0}\right\} $, it can fulfills
these conditions with many other systems of bijections which fulfill
conditions B1) and B2).

Contrariwise, if we have a system of bijections $S=\left\{
s_{B}:B\rightarrow B\mid B\in \mathrm{Ob}\Theta ^{0}\right\} $ which
fulfills conditions B1) and B2) then we can define a strongly stable
automorphism $\Phi \left( S\right) =\Phi $\ of the category $\Theta ^{0}$ by
this way: $\Phi $ preserves all objects of the category $\Theta ^{0}$ and
acts on its morphisms according formula (\ref{biject_action}) with $%
s_{B}^{\Phi }=s_{B}$. Of course, two different system of bijections $%
S_{1}=\left\{ s_{1,B}^{\Phi }:B\rightarrow B\mid B\in \mathrm{Ob}\Theta
^{0}\right\} $ and $S_{2}=\left\{ s_{2,B}^{\Phi }:B\rightarrow B\mid B\in 
\mathrm{Ob}\Theta ^{0}\right\} $ which fulfill conditions B1) and B2) can
provide by formula (\ref{biject_action}) the same action on homomorphisms
and, so, they will provide the same strongly stable automorphism of the
category $\Theta ^{0}$.

\section{Automorphic equivalence of one-sorted algebras.\label{autom_equiv}}

\setcounter{equation}{0}

In this Section we assume that there is a strongly stable automorphisms $%
\Phi $ of the category $\Theta ^{0}$. It fulfills condition A2) and A3) with
the system of bijections $S=\{s_{B}^{\Phi }:B\rightarrow B\mid B\in \mathrm{%
Ob}\Theta ^{0}\}$, which fulfills conditions B1) and B2). Then we have a
system of words $\mathfrak{W}\left( S\right) =\left\{ w_{\omega }\left(
x\right) \mid \omega \in \Omega \right\} $ which fulfills conditions Op1)
and Op2) with $\sigma _{B}=s_{B}^{\Phi }$. As above, $\omega _{H}^{\ast }$
is an operation induced on $H$ by the word $w_{\omega }\left( x\right) $ ($%
H\in \Theta $) and $H^{\ast }$ is the algebra which has same domain as the
algebra $H$ and the operations $\left\{ \omega _{H}^{\ast }\mid \omega \in
\Omega \right\} $. By \textbf{Proposition \ref{inthe_variety}}, $H^{\ast
}\in \Theta $.

\begin{proposition}
\label{clos_congr}\textbf{\hspace{-0.08in}. }\textit{Let }$B\in \mathrm{Ob}%
\Theta ^{0}$\textit{, }$H\in \Theta $\textit{\ and }$T\subset B^{2}$\textit{%
. If }$T$\textit{\ is an }$H$\textit{-closed congruence, then }$\sigma
_{B}^{-1}T$\textit{\ is an }$H^{\ast }$\textit{-closed congruence. }
\end{proposition}

\begin{proof}
We shall consider the diagram%
\[
\begin{array}{ccc}
B & \underrightarrow{\sigma _{B}} & B \\ 
\downarrow \psi &  & \varphi \downarrow \\ 
H^{\ast } &  & H%
\end{array}%
. 
\]%
If $\varphi \in \mathrm{Hom}\left( B,H\right) $ then, $\varphi \in \mathrm{%
Hom}\left( B^{\ast },H^{\ast }\right) $, $\sigma _{B}\in \mathrm{Hom}\left(
B,B^{\ast }\right) $, so $\varphi \sigma _{B}\in \mathrm{Hom}\left(
B,H^{\ast }\right) $. If $\psi \in \mathrm{Hom}\left( B,H^{\ast }\right) $
then $\psi \sigma _{B}^{-1}\in \mathrm{Hom}\left( B^{\ast },H^{\ast }\right) 
$. $\mathfrak{W}\left( S\right) $ fulfills conditions Op1) and Op2), so we
can use \textbf{Corollary \ref{invers_homomorph}} from \textbf{Proposition %
\ref{invers}} and conclude that, $\psi \sigma _{B}^{-1}\in \mathrm{Hom}%
\left( B,H\right) $.

If $T$\ is a $H$-closed congruence, then $\bigcap\limits_{\varphi \in
T_{H}^{\prime }}\ker \varphi =T$. If $\varphi \in \mathrm{Hom}\left(
B,H\right) $ and $\ker \varphi \supset T$, then $\varphi \sigma _{B}\in 
\mathrm{Hom}\left( B,H^{\ast }\right) $ and $\ker \varphi \sigma _{B}=\sigma
_{B}^{-1}\ker \varphi \supset \sigma _{B}^{-1}T$. So $\bigcap\limits_{\psi
\in \left( \sigma _{B}^{-1}T\right) _{H^{\ast }}^{\prime }}\ker \psi \subset
\bigcap\limits_{\varphi \in T_{H}^{\prime }}\ker \varphi \sigma _{B}=\sigma
_{B}^{-1}T$.
\end{proof}

\begin{remark}
\label{inv_clos_congr}\textbf{\hspace{-0.08in}. }\emph{Similar we can prove
that if }$T\subset B^{2}$\emph{\ and }$T$\emph{\ is a }$H^{\ast }$\emph{%
-closed congruence, then }$\sigma _{B}T$\emph{\ is a }$H$\emph{-closed
congruence.}
\end{remark}

\begin{theorem}
\label{H_H*}\textbf{\hspace{-0.08in}. }\textit{Automorphism }$\Phi ^{-1}$%
\textit{\ provides the automorphic equivalence of algebras }$H$\textit{\ and 
}$H^{\ast }$\textit{\ in the variety }$\Theta $\textit{.}
\end{theorem}

\begin{proof}
For every $B\in \mathrm{Ob}\Theta ^{0}$ we take the monotone bijections $%
Cl_{H}\left( B\right) \ni T\rightarrow \sigma _{B}^{-1}T\in Cl_{H^{\ast
}}\left( B\right) $. And now we need to prove that if $B_{1},B_{2}\in 
\mathrm{Ob}\Theta ^{0}$, $\mu _{1},\mu _{2}\in \mathrm{Hom}\left(
B_{1},B_{2}\right) $, $T\in Cl_{H}\left( B_{2}\right) $, $\tau
:B_{2}\rightarrow B_{2}/T$ and $\widetilde{\tau }:B_{2}\rightarrow
B_{2}/\sigma _{B_{2}}^{-1}T$ are natural homomorphisms and 
\begin{equation}
\tau \mu _{1}=\tau \mu _{2},  \label{before_functor}
\end{equation}%
then%
\begin{equation}
\widetilde{\tau }\Phi ^{-1}\left( \mu _{1}\right) =\widetilde{\tau }\Phi
^{-1}\left( \mu _{2}\right) .  \label{after_functor}
\end{equation}%
We denote $\sigma _{i}=\sigma _{B_{i}}$ ( $i=1,2$). From (\ref%
{before_functor}) we have $\left( \mu _{1}\sigma _{1}\left( b\right) ,\mu
_{2}\sigma _{1}\left( b\right) \right) \in T$ for every $b\in B$, $\left(
\sigma _{2}^{-1}\mu _{1}\sigma _{1}\left( b\right) ,\sigma _{2}^{-1}\mu
_{2}\sigma _{1}\left( b\right) \right) =\left( \Phi ^{-1}\left( \mu
_{1}\right) \left( b\right) ,\Phi ^{-1}\left( \mu _{2}\right) \left(
b\right) \right) \in \sigma _{2}^{-1}T$. Therefore (\ref{after_functor}) is
fulfilled.
\end{proof}

\section{Automorphic equivalence and geometric equivalence.\label%
{autom_geom_equiv}}

\setcounter{equation}{0}

\begin{theorem}
\label{main}\textbf{\hspace{-0.08in}. }\textit{Let the algebras }$H_{1}$%
\textit{\ and }$H_{2}$\textit{\ belongs to the variety }$\Theta $\textit{.
They are automorphically equivalent in }$\Theta $\textit{\ if and only if\
the algebra }$H_{1}$\textit{\ geometrically equivalent to the algebra }$%
H_{2}^{\ast }$\textit{, where }$H_{2}^{\ast }$\textit{\ is an algebra which
has the same domain as the algebra }$H_{2}$\textit{\ and its operations are
induced by some system of words }$\left\{ w_{\omega }\left( x\right) \mid
\omega \in \Omega \right\} $\textit{, which fulfills conditions Op1) and
Op2).}
\end{theorem}

\begin{proof}
Let an automorphism $\Psi :\Theta ^{0}\rightarrow \Theta ^{0}$ provide the
automorphic equivalence of the algebras $H_{1}$ and $H_{2}$. By \cite[%
Theorem 2]{PlotkinZhitAutCat} $\Psi $ can be decomposed: $\Psi =\Upsilon
\Phi $, where $\Upsilon $ is an inner automorphism of $\Theta ^{0}$ and $%
\Phi $ is a strongly stable one. $S=\left\{ s_{B}^{\Phi }:B\rightarrow B\mid
B\in \mathrm{Ob}\Theta ^{0}\right\} $ is a system of bijections described in
conditions A2) and A3). So the system $S$ fulfills the conditions B1) and
B2). We have the system of words $\mathfrak{W}\left( S\right) =\left\{
w_{\omega }\left( x\right) \mid \omega \in \Omega \right\} $. Let $%
H_{2}^{\ast }$\ be an algebra which has the same domain as the algebra $%
H_{2} $\ and its operations are induced by the system $\mathfrak{W}\left(
S\right) $. By \textbf{Theorem \ref{H_H*}}, $\Phi ^{-1}$ provides the
automorphic equivalence of the algebras $H_{2}$\ and $H_{2}^{\ast }$. Hence
(see Introduction), an automorphism $\Lambda =\Phi ^{-1}\Upsilon \Phi =\Phi
^{-1}\Psi $ provides the automorphic equivalence of algebras $H_{1}$\ and $%
H_{2}^{\ast }$. But, $\Lambda $ is an inner automorphism (see Introduction),
hence $H_{1}$ and $H_{2}^{\ast }$ are geometrically equivalent.

Let us have a system of words $W=\left\{ w_{\omega }\left( x\right) \mid
\omega \in \Omega \right\} $, which fulfills conditions Op1) and Op2) and $%
H_{1}$ is geometrically equivalent to $H_{2}^{\ast }$ ($H_{2}^{\ast }$ is
the algebra with the same domain as the algebra $H_{2}$ and its operations
are induced by $\left\{ w_{\omega }\left( x\right) \mid \omega \in \Omega
\right\} $). By \textbf{Proposition \ref{sys_biject}}, the system of
bijections $\mathfrak{S}\left( W\right) =\left\{ \sigma _{B}:B\rightarrow
B\mid B\in \mathrm{Ob}\Theta ^{0}\right\} $\ fulfills conditions B1) and
B2), where $\sigma _{B}$ is an isomorphism which is described in condition
Op2). We can consider the system of words $\mathfrak{W}\left( \mathfrak{S}%
\left( W\right) \right) $ and the automorphism $\Phi \left( \mathfrak{S}%
\left( W\right) \right) =\Phi $ which acts on morphisms of the category $%
\Theta ^{0}$ by bijections $\left\{ \sigma _{B}:B\rightarrow B\mid B\in 
\mathrm{Ob}\Theta ^{0}\right\} $. By \textbf{Proposition \ref{operators}} $%
\mathfrak{W}\left( \mathfrak{S}\left( W\right) \right) =W$, so, by \textbf{%
Theorem \ref{H_H*}}, the automorphism $\Phi ^{-1}$ provides the automorphic
equivalence of the algebras $H_{2}$\ and $H_{2}^{\ast }$. But automorphic
equivalence is a symmetric and transitive relation.
\end{proof}

\section{Acknowledgements.}

\setcounter{equation}{0}

I dedicate this paper to the 80th birthday of Prof. B.Plotkin. He motivated
this research. He was very heedful to it and gave my some very useful
recommendation, how to explain my results in this paper. Dr. G. Zhitomirski
was also very heedful to my research. Of course, I must appreciative to
Prof. B.Plotkin and Dr. G. Zhitomirski for theirs marvelous paper \cite%
{PlotkinZhitAutCat}. This research could not take place without theirs deep
results. Prof. S. Margolis also was very heedful to my research. Very useful
discussions with Prof. L. Rowen and Dr. E. Plotkin also help me in the
writing of this paper.

\end{document}